\documentclass[12pt]{article}

\usepackage{a4wide}

\usepackage[intlimits]{amsmath}
\usepackage{amssymb,amsthm}

\usepackage[cp1251]{inputenc}
\usepackage[english]{babel}

\newtheorem{thm}{Theorem}[section]

\newtheorem{lem}[thm]{Lemma}

\theoremstyle{definition}

\theoremstyle{remark}
\newtheorem{rem}[thm]{Remark}
\numberwithin{equation}{section}

\begin{document}\Large

\title{Asymptotic Expansion for the Functional of Markovian
Evolution in $R^d$ in the Circuit of Diffusion Approximation}

\author{I.V.Samoilenko\\
Institute of Mathematics,\\ Ukrainian National Academy of
Sciences,\\ 3 Tereshchenkivs'ka, Kyiv, 01601, Ukraine\\
isamoil@imath.kiev.ua}

\maketitle

\abstract{Is studied asymptotic expansion for solution of
singularly perturbed equation for functional of Markovian
evolution in $R^d$. The view of regular and singular parts of
solution is found. \\ {\bf Mathematics Subject Classification
(2000):} primary 60J25,secondary 35C20. \\ {\bf Keywords:} random
evolution, singularly perturbed equation, asymptotic expansion,
diffusion approximation, estimate of the remainder}

\section{\bf Introduction}
The problems of asymptotic expansion for solutions of PDE and PDE
systems were studied by many authors. A lot of references could be
found in {\rm \cite{Mark}}. There are studied, as a rule, border
problems and the small parameter is at the higher derivative by
$t$.

For example, in {\rm \cite{VasBut}}(p. 155) is studied the system
of first order equations with the small parameter by $t$ and $x$
that corresponds the telegraph equation.

In this work we study asymptotic expansion for solution of
singularly perturbed equation for functional of Markovian
evolution in $R^d$.

Let $x\in R^d, \xi(s)$ - an ergodic Markovian process in the set
$E=\{1,\ldots,N\}$ with the intensity matrix $Q=\{q_{ij},
i,j=\overline{1,N}\}.$

The probability of being in the $i$ - th state longer then $t$ is
$P\{\theta_i>t\}=e^{-q_it},$ where $q_i=\sum\limits_{j\neq i }
q_{ij}.$

Let $a(i)=(a_1(i),\ldots,a_d(i))$ - vector-function on $E$. We
regard a vector-function as a corresponding vector-column.

Put matrix $A=\{a_k(i), k=\overline{1,d}, i=\overline{1,N}\}.$

Let us study evolution
$$x^{\varepsilon}(t)=x+\varepsilon^{-1}\int_0^ta(\xi(s/\varepsilon^2))ds=$$ $$x+\varepsilon\int_0^{t/\varepsilon^2}
a(\xi(s))ds.$$

It's well-known {\rm \cite{Pin}}, that the functionals of
evolution, determined
 by a test-function $\overline{f}(x)\in C^{\infty}(R^d)$ (here $\overline{f}(x)$
is integrable on $R^d$ and has equal components
$\overline{f}(x)=(f(x),\ldots,f(x))$):
$u_i^{\varepsilon}(x,t)=E_i\overline{f}(x^{\varepsilon}(t)),
i=\overline{1,N}$ (here $i$ - is a start state of $\xi(s)$)
satisfy the system of Kolmogorov backward differential equations:
$$\frac{\partial}{\partial
t}{u}^{\varepsilon}({x},t)=\varepsilon^{-2}
Q{u}^{\varepsilon}({x},t)+\varepsilon^{-1}A\nabla{u}^{\varepsilon}({x},t),
\eqno(1)$$ where
${u}^{\varepsilon}({x},t)=(u_{1}^{\varepsilon}({x},t),\ldots,u_{N}^{\varepsilon}({x},t)),$
$A\nabla=diag[(a(i),\nabla),i=\overline{1,N}],$
$\nabla=(\frac{\partial}{\partial
x_1},\ldots,\frac{\partial}{\partial x_d}).$

As an example we'll describe a well-known model, where an equation
of type (1) appears.

{\bf Example 1.1:} {\it In the works {\rm \cite{Pin,Sam}}
functionals of the view
$$u_i(x,t)=E_i\overline{f}(x+v\int_0^t\overline{\tau}_{\xi(s)}ds),
i=\overline{0,n}$$ were studied. Here $\xi(u)$ - Poisson process
with parameter $\lambda$, $\xi(0)=0$, $v$ - velocity of particle's
motion, $\overline{\tau}_i, i=\overline{0,n}$ - vectors that
determine the directions of motion. The systems of Kolmogorov
backward differential equations were received for the functionals
$u_i(x,t), i=\overline{0,n}$ in case of cyclic and uniform change
of motion directions.

In a matrix form we have: $$\frac{\partial}{\partial
t}{u}^{\varepsilon}(x,t)=[\lambda
Q+vA\nabla]{u}^{\varepsilon}(x,t), \eqno(2)$$ where
${u}^{\varepsilon}(x,t)=(u_{0}^{\varepsilon}(x,t),\ldots,u_{n}^{\varepsilon}(x,t)),$
$A\nabla=diag[(\overline{\tau}_i,\nabla), i=\overline{0,n}],$
$Q=[q_{ij}, i,j=\overline{0,n}].$ Here $q_{ii}=-1,$ $q_{ii+1}=1,$
$q_{ij}=0,j\neq i,j\neq i+1$ in case of cyclic change of
directions, and $q_{ii}=-1, q_{ij}=1/n, i\neq j$ in case of
uniform change.

If we put in (2) $v=\varepsilon^{-1}, \lambda=\varepsilon^{-2}$,
where $\varepsilon$ - is a small parameter, we'll have a
singularly perturbed equation of type (1):
$$\frac{\partial}{\partial
t}{u}^{\varepsilon}(x,t)=[\varepsilon^{-2}Q+
\varepsilon^{-1}A\nabla]{u}^{\varepsilon}(x,t).$$

Initial condition
${u}^{\varepsilon}(x,0)=\overline{f}(x):=(f(x),\ldots,f(x)).$}

Equations of type (1) were also studied in the works {\rm
\cite{Kor,KorPenTur}}, partially in {\rm \cite{KorPenTur}} for the
distribution of absorption time of Markov chain with continuous
time that depends on small parameter $\varepsilon$ the following
equation was received
$\varepsilon\frac{d}{dx}u^{\varepsilon}(x)=(Q-\varepsilon
G)u^{\varepsilon}(x), Q=P-I.$ Asymptotic expansion of its solution
was found there.

In this work we study system (1) with the second order
singularity. This problem has interesting probabilistic sense:
hyperbolic equation of high degree, corresponding system (2) (see
{\rm \cite{Sam}}) becomes parabolic equation of Wiener process in
hydrodynamic limit, when $\varepsilon\to 0$. The fact that
solution of (2) in hydrodynamic limit tends to the functional of
Wiener process is well-known and studied, for example, in {\rm
\cite{KorTur}}.

To find asymptotic expansion of the solution of (1) we use the
method proposed in {\rm \cite{VasBut}}. The solution consists two
parts - regular terms and singular terms - which are determined by
different equations. Asymptotic expansion lets not only determine
the terms of asymptotic, but to see the velocity of convergence in
hydrodynamic limit.

Besides, when studying this problem, we improved the algorithm of
asymptotic expansion. Partially, the initial conditions for the
regular terms of asymptotic are determined without the use of
singular terms, i.e. the regular part of the solution may be found
by a separate recursive algorithm; scalar part of the regular term
is found and without the use of singular terms. These and other
improves of the algorithm are pointed later.

\section{\bf Asymptotic expansion of the solution}
Let $P(t)=e^{Qt}=\{p_{ij}(t); i,j=\overline{1,N}\}.$

Put $\pi_j=\lim\limits_{t\to\infty}p_{ij}(t)$ and
$-R_0=\{\int_0^{\infty}(p_{ij}(t)-\pi_j)dt;
i,j=\overline{1,N}\}=\{r_{ij}; i,j=\overline{1,N}\}.$

Let $\Pi$ be a projecting matrix on the null-space $N_Q$ of the
matrix $Q$. For any vector $g$ we have $\Pi
g=\widehat{g}\mathbf{1},$ where
$\widehat{g}=\sum\limits_{i=1}^Ng_i\pi_i,
\mathbf{1}=(1,\ldots,1).$ Then for the matrix $Q$ the following
correlation is true $\Pi Q\Pi=0$ (see {\rm \cite{KorTur}}, chapter
3).

Let the matrix $A$ satisfy balance condition: $$\Pi A\Pi=0.$$

We put: $$R_0A\nabla=\{r_{ij}(a(j),\nabla), i,j=\overline{1,N}
\}=\left\{\sum_{k=1}^d r_{ij}a_k(j)\frac{\partial}{\partial x_k} ,
i,j=\overline{1,N}\right\},$$ $$A\nabla R_0=\{(a(i),\nabla)r_{ij},
i,j=\overline{1,N} \}=\left\{\sum_{k=1}^d
a_k(i)r_{ij}\frac{\partial}{\partial x_k} ,
i,j=\overline{1,N}\right\},$$ $$A\nabla
R_0A\nabla=\{(a(i),\nabla)r_{ij}(a(j),\nabla), i,j=\overline{1,N}
\}=$$ $$\left\{\sum_{k=1}^d\sum_{l=1}^d
a_k(i)r_{ij}a_l(j)\frac{\partial}{\partial
x_k}\frac{\partial}{\partial x_l} , i,j=\overline{1,N}\right\},$$
$$exp_0(Qt):=e^{Qt}-\Pi,$$
$$\widehat{a}_{kl}=\sum_{i,j=1}^N\pi_ia_k(i)r_{ij}a_l(j)\pi_j, $$
here, following {\rm \cite{Pin}}, we need the condition:
$$\widehat{a}_{kl}>0.$$

\begin{thm} The solution of equation (1) with initial condition
${u}^{\varepsilon}(x,0)=\overline{f}(x)$, where
$\overline{f}(x)\in C^{\infty}(R^d)$ and integrable on $R^d$ has
asymptotic expansion
$${u}^{\varepsilon}(x,t)={u}^{(0)}(x,t)+\sum_{n=1}^{\infty}\varepsilon^n\left({u}^{(n)}(x,t)+
{v}^{(n)}\left(x,t/\varepsilon^2\right)\right).\eqno(3)$$

Regular terms of the expansion are: ${u}^{(0)}(x,t)$ - the
solution of equation
$$\frac{\partial}{\partial
t}{u}^{(0)}(x,t)=\sum_{k,l=1}^d\widehat{a}_{kl}\frac{\partial^2
{u}^{(0)}(x,t)}{\partial x_k \partial x_l} \eqno(4)$$ with initial
condition ${u}^{(0)}(x,0)=\overline{f}(x),$
$${u}^{(1)}(x,t)=R_0A\nabla{u}^{(0)}(x,t)=\left[\sum_{k=1}^d\sum_{j=1}^N
r_{ij}a_k(j)\frac{\partial u^{(0)}_{j}(x,t)}{\partial x_k} ,
i=\overline{1,N} \right]$$ for $k\geq 2:$
$${u}^{(k)}(x,t)=R_0\left[\frac{\partial}{\partial
t}{u}^{(k-2)}(x,t)-A\nabla{u}^{(k-1)}(x,t)\right]+c^{(k)}(t):=$$
$$ :=R_0\Phi
\left[{u}^{(k-2)}(x,t),{u}^{(k-1)}(x,t)\right]+c^{(k)}(t),$$ where
$$c^{(k)}(t)\in N_{Q}, c^{(k)}(t)=c^{(k)}(0)+\int_0^t\widehat{L_k}c^{(0)}(s)ds,$$
here
$$c^{(0)}(t)={u}^{(0)}(x,t),
L_0=\left\{\sum_{k=1}^d\sum_{l=1}^d
a_k(i)r_{ij}a_l(j)\frac{\partial}{\partial
x_k}\frac{\partial}{\partial x_l} , i,j=\overline{1,N} \right\}$$
$$ \widehat{L_k}=\Pi L_k\Pi, L_k=(-1)^{k+1}(A\nabla
R_0)^k\pounds_0, k\geq 1, $$ $$
\pounds_0=\left\{\frac{\partial}{\partial
t}-\sum_{k=1}^d\sum_{l=1}^d
a_k(i)r_{ij}a_l(j)\frac{\partial}{\partial
x_k}\frac{\partial}{\partial x_l} , i,j=\overline{1,N} \right\}.$$

The singular terms of the expansion have the
view:$${v}^{(1)}(x,t)=exp_0({Qt})A\nabla \overline{f}(x),$$ for
$k>1:$
$${v}^{(k)}(x,t)=exp_0(Qt){v}^{(k)}(x,0)+\int_0^texp_0(Q(t-s))A\nabla{v}^{(k-1)}(x,s)ds-$$ $$\Pi
\int_t^{\infty}A\nabla{v}^{(k-1)}(x,s)ds.$$ Initial conditions:
$$c^{(0)}(0)=\overline{f}(x),$$
$$u^{(1)}(x,0)=R_0A\nabla \overline{f}(x),
{v}^{(1)}(x,0)=-\frac{1}{2}A\nabla\Pi \overline{f}(x),$$ for
$k>1:$
$${v}^{(k)}(x,0)=\Phi\left[{u}^{(k-2)}(x,0),{u}^{(k-1)}(x,0)\right],$$
$$c^{(k)}(0)=-A\nabla\widetilde{v}^{(k-1)}(x,0),$$
where $ \widetilde{v}^{(1)}(x,0)=-R_0A\nabla \overline{f}(x),$
$$\widetilde{v}^{(k)}(x,0)=R_0\Phi
\left[{u}^{(k-2)}(x,0),{u}^{(k-1)}(x,0)\right]+R_0A\nabla\widetilde{v}^{(k-1)}(x,0)+$$
$$\Pi
A\nabla(\widetilde{v}^{(k-1)}(x,\lambda))'_{\lambda}|_{\lambda=0},$$
$$(\widetilde{v}^{(k)}(x,\lambda))'_{\lambda}|_{\lambda=0}=R^2_0\Phi
\left[{u}^{(k-2)}(x,0),{u}^{(k-1)}(x,0)\right]+R_0^2Q_{1}\widetilde{v}^{(k-1)}(x,0)+$$
$$
R_0A\nabla(\widetilde{v}^{(k-1)}(x,\lambda))'_{\lambda}|_{\lambda=0}.$$

\end{thm}

\begin{rem} The initial conditions for the regular terms of
asymptotic are determined without the use of singular terms, i.e.
the regular part of the solution may be found by a separate
recursive algorithm (comp. with {\rm \cite{KorPenTur}}).
\end{rem}

\begin{rem} In case of evolution described in Example 1
equation (4) has the view:
$$\frac{\partial}{\partial
t}{u}^{(0)}(x,t)=\frac{1}{(n+1)^2}\triangle{u}^{(0)}(x,t)$$ with
initial condition ${u}^{(0)}(x,0)=\overline{f}(x).$

Solution of this problem in the class of integrable and infinitely
differentiable functions of exponential growth is:
$${u}^{(0)}({x},t)=(2\pi
t)^{-\frac{n}{2}}\frac{1}{(n+1)^2}\int_{R^n}e^{-(n+1)^2\frac{<({x}-{y}),({x}-{y})>}{2t}}
\overline{f}({y})d{y}.$$
\end{rem}

{\it Proof of Theorem 2.1:} Let us substitute the solution
${u}^{\varepsilon}(x,t)$ in the view (3) to the equation (1) and
equal the terms at $\varepsilon$ degrees. We'll have the system
for the regular terms of asymptotic:$$
  \begin{cases}
    Q{u}^{(0)}=0 \\
    Q{u}^{(1)}+A\nabla {u}^0=0 \\
    Q{u}^{(k)}=\frac{\partial}{\partial
    t}{u}^{(k-2)}-A\nabla{u}^{(k-1)}, k\geq 2
  \end{cases}\eqno(5)
$$   and for the singular terms: $$
  \begin{cases}
    \frac{\partial}{\partial
    t}{v}^{(1)}=Q{v}^{(1)} \\
    \frac{\partial}{\partial
    t}{v}^{(k)}-Q{v}^{(k)}=A\nabla{v}^{(k-1)},
    k>1.
  \end{cases} \eqno(6)
$$

From (5) we have: ${u}^{(0)}\in N_{Q},
{u}^{(1)}=R_0A\nabla{u}^{(0)}+c^{(1)}(t).$ For ${u}^{(2)}$ we
receive: $Q{u}^{(2)}=\frac{\partial}{\partial
t}{u}^{(0)}-A\nabla{u}^{(1)}=\frac{\partial}{\partial
t}{u}^{(0)}-A\nabla R_0A\nabla{u}^{(0)}=\frac{\partial}{\partial
t}{u}^{(0)}-L_0{u}^{(0)}.$

The solvability condition for ${u}^{(2)}$ has the view: $$\Pi Q\Pi
{u}^{(2)}=0= \frac{\partial}{\partial t}{u}^{(0)}-\Pi L_0\Pi
{u}^{(0)}.$$

So, we have equation (4) for ${u}^{(0)}(x,t)$.

We note that in {\rm \cite{KorPenTur}} solvability condition is
written for the equation that contains the terms ${u}^{(0)}(x,t)$
and ${u}^{(1)}(x,t)$. In this work we have to express
${u}^{(1)}(x,t)$ through ${u}^{(0)}(x,t)$ and only then to write
down solvability condition for the equation that contains the
terms ${u}^{(0)}(x,t)$ and ${u}^{(2)}(x,t)$.

For ${u}^{(1)}$ we have:
$${u}^{(1)}=R_0A\nabla{u}^{(0)}+c^{(1)}(t).$$

Using the last equation from (5) we receive:
$${u}^{(k)}(x,t)=R_0\left[\frac{\partial}{\partial
t}{u}^{(k-2)}(x,t)-A\nabla{u}^{(k-1)}(x,t)\right]+c^{(k)}(t):=$$
$$ :=R_0\Phi
\left[{u}^{(k-2)}(x,t),{u}^{(k-1)}(x,t)\right]+c^{(k)}(t),$$ where
$c^{(k)}(t)\in N_{Q}.$

To find $c^{(k)}(t)$ we'll use the fact that ${u}^{(0)}\in N_{Q}.$
Let us put $c^{(0)}(t)={u}^{(0)}(x,t)$. From the equation
$Q{u}^{(2)}=\frac{\partial}{\partial
t}c^{(0)}(t)-L_0c^{(0)}(t)=\pounds_0c^{(0)}(t)$ we have
$${u}^{(2)}=R_0\pounds_0c^{(0)}(t).$$

For ${u}^{(3)}$: $$Q{u}^{(3)}=\frac{\partial}{\partial
t}c^{(1)}(t)-A\nabla{u}^{(2)}=(c^{(1)}(t))'-A\nabla
R_0\pounds_0c^{(0)}(t)=\pounds_1c^{(0)}(t).$$ From the solvability
condition $\Pi Q\Pi {u}^{(3)}=0=\frac{\partial}{\partial
t}c^{(1)}(t)-\Pi A\nabla R_0 \pounds_0 \Pi
c^{(0)}(t)=(c^{(1)}(t))'-\widehat{L}_1c^{(0)}(t)$ we find:
$$c^{(1)}(t)=c^{(1)}(0)+\int_0^t\widehat{L}_1c^{(0)}(s)ds,$$  and
${u}^{(3)}=R_0\pounds_1c^{(0)}(t),$ where
$\pounds_1=(-L_1)c^{(0)}(t),$ as soon as $R_0\widehat{L}_1=0.$

By induction:
$$c^{(k)}(t)=c^{(k)}(0)+\int_0^t\widehat{L}_kc^{(0)}(s)ds,$$ where
$\widehat{L}_k=\Pi L_k\Pi , L_k=(-1)^{k+1}(A\nabla R_0
)^k\pounds_0, \pounds_0=\frac{\partial}{\partial t}-L_0, k\geq 2.$

In contrast to {\rm \cite{KorPenTur}}, where the equations for
$c^{(k)}(t)$ were found, in this work we may find $c^{(k)}(t)$
explicitly through $c^{(0)}(t)$.

For the singular terms we have from (6):
$${v}^{(1)}(x,t)=exp_0(Qt){v}^{(1)}(x,0).$$

Here we should note that the ordinary solution
${v}^{(1)}(x,t)=exp(Qt){v}^{(1)}(x,0)$ is corrected by the term
$-\Pi{v}^{(1)}(x,0)$ in order
 to receive the following limit
 $\lim\limits_{t\to
 \infty}{v}^{(1)}(x,t)=0$.
This limit is true for all singular terms due to uniform
ergodicity of switching
 Markovian process.

For the homogenous part of the second equation of the system we
have the following solution:
$${v}^{(k)}(x,t)=exp_0(Qt){v}^{(k)}(x,0).$$

But as soon as the equation is not homogenous the corresponding
solution should be
$${v}^{(k)}(x,t)=exp_0(Qt){v}^{(k)}(x,0)+
\int_0^texp_0(Q(t-s))A\nabla{v}^{(k-1)}(x,s)ds.$$

But here we should again correct the solution, in order
 to receive the limit
 $\lim\limits_{t\to
 \infty}{v}^{(k)}(x,t)=0$, by the term $-
\Pi\int_t^{\infty}A\nabla {v}^{(k-1)}(x,s)ds.$

And so the solution is:
$${v}^{(k)}(x,t)=exp_0(Qt){v}^{(k)}(x,0)+
\int_0^texp_0(Q(t-s))A\nabla{v}^{(k-1)}(x,s)ds-$$
$$\Pi\int_t^{\infty}A\nabla {v}^{(k-1)}(x,s)ds.$$

We should finally find the initial conditions for the regular and
singular terms.

We put $c^{(0)}(t)={u}^{(0)}(x,t)$, so
$c^{(0)}(0)={u}^{(0)}(x,0)=\overline{f}(x).$

From the initial condition for the solution
${u}^{\varepsilon}(x,0)={u}^{(0)}(x,0)=(f(x),\ldots,f(x))$, we
have ${u}^{(k)}(x,0)+{v}^{(k)}(x,0)=0, k\geq 1.$ Let us rewrite
this equation for the null-space $N_Q$ of matrix
$Q$:$$\Pi{u}^{(k)}(x,0)+\Pi{v}^{(k)}(x,0)=0, k\geq 1, \eqno(7)$$
and the space of values $R_Q$:
$$(I-\Pi){u}^{(k)}(x,0)+(I-\Pi){v}^{(k)}(x,0)=0, k\geq 1. \eqno(8)$$

As we proved for $k>1$:
$${u}^{(k)}(x,0)=R_0\Phi[{u}^{(k-2)}(x,0),{u}^{(k-1)}(x,0)]+c^{(k)}(0)=$$ $$=(I-\Pi)\Phi[{u}^{(k-2)}(x,0),{u}^{(k-1)}(x,0)]+\Pi c^{(k)}(0),$$
$${v}^{(k)}(x,0)=(I-\Pi){v}^{(k)}(x,0)-\Pi\int_0^{\infty}A\nabla{v}^{(k-1)}(x,s)ds.$$

Functions ${v}^{(k-1)}(x,s),{u}^{(k-2)}(x,0),{u}^{(k-1)}(x,0)$ are
known from the previous steps of induction. So, we've found
$\Pi{v}^{(k)}(x,0)$ in (7) and $(I-\Pi){u}^{(k)}(x,0)$ in (8).

Now we may use the correlations (7), (8) to find the unknown
initial conditions:
$$c^{(k)}(0)=-\int_0^{\infty}A\nabla{v}^{(k-1)}(x,s)ds,$$
$${v}^{(k)}(x,0)=\Phi[{u}^{(k-2)}(x,0),{u}^{(k-1)}(x,0)].$$

In {\rm \cite{KorPenTur}} an analogical correlation was found for
$c^{(k)}(0)$. To find $c^{(k)}(0)$ explicitly and without the use
of singular terms we'll find Laplace transform for the singular
term. The following lemma is true.

\begin{lem} Laplace transform for the singular term of asymptotic
expansion
$$\widetilde{v}^{(k)}(x,\lambda)=\int_0^{\infty}e^{-\lambda s}
{v}^{(k)}(x,s)ds$$ has the view:

$$\widetilde{v}^{(1)}(x,\lambda)=(\lambda-\Pi+(R_0+\Pi)^{-1})^{-1}[-R_0A\nabla
\overline{f}(x)],$$

$$\widetilde{v}^{(k)}(x,\lambda)=(\lambda-\Pi+(R_0+\Pi)^{-1})^{-1}\Phi
\left[{u}^{(k-2)}(x,0),{u}^{(k-1)}(x,0)\right]+$$
$$(\lambda-\Pi+(R_0+\Pi)^{-1})^{-1}A\nabla
\widetilde{v}^{(k-1)}(x,\lambda)+\frac{1}{\lambda}\Pi
A\nabla[\widetilde{v}^{(k-1)}(x,\lambda)-\widetilde{v}^{(k-1)}(x,0)],$$
where $$\widetilde{v}^{(1)}(x,0)=-R_0A\nabla \overline{f}(x),$$

$$(\widetilde{v}^{(1)}(x,\lambda))'_{\lambda}|_{\lambda=0}=-R_0^2A\nabla\Pi
\overline{f}(x),$$

$$\widetilde{v}^{(k)}(x,0)=R_0\Phi
\left[{u}^{(k-2)}(x,0),{u}^{(k-1)}(x,0)\right]+R_0A\nabla\widetilde{v}^{(k-1)}(x,0)+$$
$$\Pi
A\nabla(\widetilde{v}^{(k-1)}(x,\lambda))'_{\lambda}|_{\lambda=0},$$

$$(\widetilde{v}^{(k)}(x,\lambda))'_{\lambda}|_{\lambda=0}=R^2_0\Phi
\left[{u}^{(k-2)}(x,0),{u}^{(k-1)}(x,0)\right]+R_0^2Q_{1}\widetilde{v}^{(k-1)}(x,0)+$$
$$
R_0A\nabla(\widetilde{v}^{(k-1)}(x,\lambda))'_{\lambda}|_{\lambda=0}.$$
\end{lem}

{\it Proof.}
$$\widetilde{v}^{(1)}(x,\lambda)=\int_0^{\infty}e^{-\lambda s}
{v}^{(1)}(x,s)ds=\int_0^{\infty}e^{-\lambda s} [e^{Qs}-\Pi]ds
v^{(1)}(x,0)=$$ $$=(\lambda-\Pi+(R_0+\Pi)^{-1})^{-1}[-A\nabla
\overline{f}(x)],$$ where the correlation for the resolvent was
found in {\rm \cite{KorTur}}.

$$\widetilde{v}^{(1)}(x,0)=-R_0A\nabla\overline{f}(x),$$
$$(\widetilde{v}^{(1)}(x,\lambda))'_{\lambda}|_{\lambda=0}=\lim\limits_{\lambda\to 0}\frac{R(\lambda)-R_0}{\lambda}[-A\nabla
\overline{f}(x)]=-R_0^2A\nabla \overline{f}(x).$$

For the next terms we have:
$$\widetilde{v}^{(k)}(x,\lambda)=(\lambda-\Pi+(R_0+\Pi)^{-1})^{-1}\Phi
\left[{u}^{(k-2)}(x,0),{u}^{(k-1)}(x,0)\right]+$$
$$(\lambda-\Pi+(R_0+\Pi)^{-1})^{-1}A\nabla
\widetilde{v}^{(k-1)}(x,\lambda)+\frac{1}{\lambda}\Pi
A\nabla[\widetilde{v}^{(k-1)}(x,\lambda)-\widetilde{v}^{(k-1)}(x,0)],$$
here the last term was found using the following correlation:
$$\int_0^{\infty}e^{-\lambda s }\int_s^{\infty}
A\nabla{v}^{(k-1)}(x,\tau)d\tau
ds=\int_0^{\infty}\int_0^{\tau}e^{-\lambda s
}A\nabla{v}^{(k-1)}(x,\tau)ds d\tau=$$
$$\int_0^{\infty}\left(-\frac{1}{\lambda}\right)(e^{-\lambda\tau}
-1)A\nabla{v}^{(k-1)}(x,\tau)d\tau =\frac{1}{\lambda}
A\nabla[\widetilde{v}^{(k-1)}(x,\lambda)-\widetilde{v}^{(k-1)}(x,0)].$$

So, $$\widetilde{v}^{(k)}(x,0)=R_0\Phi
\left[{u}^{(k-2)}(x,0),{u}^{(k-1)}(x,0)\right]+R_0A\nabla\widetilde{v}^{(k-1)}(x,0)+\Pi
A\nabla(\widetilde{v}^{(k-1)}(x,\lambda))'_{\lambda}|_{\lambda=0},$$

$$(\widetilde{v}^{(k)}(x,\lambda))'_{\lambda}|_{\lambda=0}=R^2_0\Phi
\left[{u}^{(k-2)}(x,0),{u}^{(k-1)}(x,0)\right]+R_0^2Q_{1}\widetilde{v}^{(k-1)}(x,0)+$$
$$R_0A\nabla(\widetilde{v}^{(k-1)}(x,\lambda))'_{\lambda}|_{\lambda=0}-\lim\limits_{\lambda\to 0}\left\{\frac{1}{\lambda^2}\Pi
A\nabla[\widetilde{v}^{(k-1)}(x,\lambda)-\widetilde{v}^{(k-1)}(x,0)]-\right.$$
$$\left.\frac{1}{\lambda}\Pi
A\nabla(\widetilde{v}^{(k-1)}(x,\lambda))'_{\lambda}\right\},$$
where the last limit tends to 0.

Lemma is proved.

So, the obvious view of the initial condition for the $c^{(k)}(t)$
is:
$$c^{(k)}(0)=-A\nabla\widetilde{v}^{(k-1)}(x,0).$$

Theorem is proved.

\section{\bf Estimate of the remainder}
Let function $f(x,i)$ in the definition of the functional
${u}^{\varepsilon}(x,t)$ belongs to Banach space of twice
continuously differentiable by $x$ functions $C^2(R^d\times E)$.

Let us write (1) in the view
$$\tilde{u}^{\varepsilon}_2(x,t)={u}^{\varepsilon}(x,t)-{u}^{\varepsilon}_2(x,t)\eqno(9)$$
where ${u}^{\varepsilon}_2(x,t)=
{u}^{(0)}(x,t)+\varepsilon({u}^{(1)}(x,t)+{v}^{(1)}(x,t))+\varepsilon^2({u}^{(2)}(x,t)+
{v}^{(2)}(x,t)),$ and the explicit view of the functions
${u}^{(i)}(x,t), {v}^{(j)}(x,t), i=\overline{0,2}, j=1,2$ is given
in Theorem 2.1.

By theorem from {\rm \cite{KorTur}} in Banach space $C^2(R^d\times
E)$ for the generator of Markovian evolution
$L^{\varepsilon}=\varepsilon^{-2}Q+\varepsilon^{-1}A\nabla,$
exists bounded inverse operator
$(L^{\varepsilon})^{-1}=\varepsilon^2[Q+\varepsilon
A\nabla]^{-1}.$

Let us substitute the function (9) into equation (1):
$$\frac{d}{dt}\tilde{u}^{\varepsilon}-L^{\varepsilon}\tilde{u}^{\varepsilon}=
\frac{d}{dt}{u}^{\varepsilon}_2-L^{\varepsilon}{u}^{\varepsilon}_2:=\varepsilon
w^{\varepsilon}. \eqno(10)$$

Here $\varepsilon
w^{\varepsilon}=\varepsilon[\frac{d}{dt}(({u}^{(1)}+{v}^{(1)})+\varepsilon
({u}^{(2)}+{v}^{(2)}))-(\varepsilon^{-1}Q({u}^{(1)}+{v}^{(1)})+Q({u}^{(2)}+
{v}^{(2)})+A\nabla({u}^{(1)}+{v}^{(1)})+\varepsilon
A\nabla({u}^{(2)}+{v}^{(2)}))].$

The initial condition has the order $\varepsilon$, so we may write
it in the view:
$$\tilde{u}^{\varepsilon}(0)=\varepsilon
\tilde{u}^{\varepsilon}(0).$$

Let
$L_t^{\varepsilon}f(x,i)=E[f(x^{\varepsilon}(t),\xi^{\varepsilon}(t/\varepsilon^2))|x^{\varepsilon}(0)=x,
\xi^{\varepsilon}(0)=i]$ be the semigroup corresponding to the
operator $L^{\varepsilon}.$

\begin{thm} The following estimate is true for the remainder (9) of the
solution of equation (1):
$$||\tilde{u}^{\varepsilon}(t)||\leq
 \varepsilon
||\tilde{u}^{\varepsilon}(0)|| \exp\{\varepsilon
L||w^{\varepsilon}||\}, $$ where $L=2||(L^{\varepsilon})^{-1}||.$
\end{thm}

{\it Proof:} The solution of equation (10) is:
$$\tilde{u}_2^{\varepsilon}(t)=\varepsilon[L_t^{\varepsilon}
\tilde{u}^{\varepsilon}(0)+ \int_0^t
L_{t-s}^{\varepsilon}w^{\varepsilon}(s)ds].$$

For the semigroup we have
$L_t^{\varepsilon}=I+L^{\varepsilon}\int_0^tL_s^{\varepsilon}ds,$
so
$\int_0^tL_s^{\varepsilon}ds=(L^{\varepsilon})^{-1}(L_t^{\varepsilon}-I).$

Using Gronwell-Bellman inequality {\rm \cite{BainSim}}, we receive
$$||\tilde{u}^{\varepsilon}(t)||\leq \varepsilon L_t^{\varepsilon}
||\tilde{u}^{\varepsilon}(0)||\exp\{\varepsilon\int_0^t
L_s^{\varepsilon}w^{\varepsilon}(t-s)ds\}\leq \varepsilon
L_t^{\varepsilon} ||\tilde{u}^{\varepsilon}(0)|| \exp\{\varepsilon
L||w^{\varepsilon}||\}, $$ where $L=2||(L^{\varepsilon})^{-1}||.$

Theorem is proved.

\begin{rem} For the remainder of asymptotic expansion (3) of the view
$$\tilde{u}^{\varepsilon}_{N+1}(x,t):={u}^{\varepsilon}(x,t)-{u}^{\varepsilon}_{N+1}(x,t),$$
where ${u}^{\varepsilon}_{N+1}(x,t)=
{u}^{(0)}(x,t)+\sum_{n=1}^{N+1}\varepsilon^n({u}^{(n)}(x,t)+{v}^{(n)}(x,t))$
we have analogical estimate:
$$||\tilde{u}^{\varepsilon}_{N+1}(t)||\leq
 \varepsilon^N
||\tilde{u}^{\varepsilon}(0)|| \exp\{\varepsilon^N
L||w^{\varepsilon}_N||\},$$ where
$\frac{d}{dt}u^{\varepsilon}_{N+1}-L^{\varepsilon}u^{\varepsilon}_{N+1}:=\varepsilon^N
w^{\varepsilon}_N.$
\end{rem}

{\it Acknowledgements.}The author thank the Institute of Applied
Mathematics, University of Bonn for the hospitality. The financial
support by DFG  project 436 UKR 113/70/0-1  is gratefully
acknowledged.\par


\begin{thebibliography}{99}

\bibitem{BainSim} Bainov D., Simeonov P. {\it Integral inequalities
and applications,} Kluver Acad. Publ., Dordrecht, (1992), 316p.

\bibitem{Kor} Korolyuk V.S., {\it Boundary layer in asymptotic
analysis for random walks,} Theory of Stochastic Processes
\underline{1-2}, 25-36 (1998).

\bibitem{KorPenTur} Koroljuk V.S., Penev I.P., Turbin A.F., {\it Asymptotic
expansion for the distribution of absorption time of Markov
chain,} Cybernetics \underline{4}, 133-135 (1973), (in Russian).

\bibitem{KorTur} Koroljuk V.S., Turbin A.F.
{\it Mathematical foundation of  state lumping of large systems,}
Kluver Acad. Press, Amsterdam, (1990), 280p.

\bibitem{Mark} Markush I.I. {\it Development of asymptotic
methods in the theory of differential equations,} Uzhgorod,
(1975), 224 p. (in Ukrainian).

\bibitem{Pin} Pinsky M. {\it Lectures on random evolutions,} World Scientific,
Singapore, (1991), 136 p.

\bibitem{Sam} Samoilenko I.V., {\it Markovian random evolution in
$R^n$,} Rand. Operat. and Stoc. Equat. \underline{2}, 139-160,
(2001).

\bibitem{VasBut} Vasiljeva A.B., Butuzov V.F. {\it Asymptotic
methods in the theory of singular perturbations,} Vyschaja shkola,
Moscow, (1990), 208 p. (in Russian).


\end{thebibliography}
\end{document}